\documentclass[12pt]{amsart}
\usepackage{amsfonts,amsmath,amssymb,amsthm}
\usepackage{latexsym}
\usepackage{mathrsfs}
\usepackage{tikz}
\usepackage{color}

\newtheorem{prop}{Proposition}[section]

\newtheorem{corollary}{Corollary}[section]

\newtheorem{example}{Example}[section]

\newcommand{\osp}{\mathfrak{osp}}

\newcommand{\id}{\text{id}}
\newcommand{\cR}{\mathcal{R}}

\newcommand{\C}{\mathcal{C}}

\begin{document}

\title{Bannai--Ito algebras and the universal $R$-matrix of $\mathfrak{osp}(1|2)$}

\author[N.Cramp\'e]{Nicolas Cramp\'e$^{\dagger}$}
\address{$^\dagger$ Institut Denis-Poisson CNRS/UMR 7013 - Universit\'e de Tours - Universit\'e d'Orl\'eans, 
Parc de Grammont, 37200 Tours, France.}
\email{crampe1977@gmail.com}

\author[L.Vinet]{Luc Vinet$^{*}$}
\address{$^*$ Centre de recherches math\'ematiques, Universit\'e de Montr\'eal,
P.O. Box 6128, Centre-ville Station,
Montr\'eal (Qu\'ebec), H3C 3J7, Canada.}

\author[M.Zaimi]{Meri Zaimi$^{*}$}

\begin{abstract}
The Bannai-Ito algebra $BI(n)$ is viewed as the centralizer of the action of $\mathfrak{osp}(1|2)$ in the $n$-fold tensor product of the universal algebra of this Lie superalgebra. 
The generators of this centralizer are constructed with the help of the universal $R$-matrix of  $\mathfrak{osp}(1|2)$. The specific structure of the  $\mathfrak{osp}(1|2)$ embeddings to 
which the centralizing elements are attached as Casimir elements is explained. With the generators defined, the structure relations of $BI(n)$ are derived from those of $BI(3)$ by repeated action of the coproduct and using properties of the $R$-matrix and of the generators of the symmetric group  $\mathfrak S_n$.

\end{abstract}

\maketitle

%

%


\section{Introduction}
This letter explains the essential role that the universal $R$-matrix of  $\mathfrak{osp}(1|2)$ plays in the algebraic underpinnings of the Bannai--Ito algebra $BI(3)$ and its higher rank generalization $BI(n)$.

The universal Bannai--Ito algebra  $BI(3)$  is generated by the central elements $\C_1$, $\C_2$, $\C_3$ and $\C_{123}$ and three generators $\C_{12}$, $\C_{23}$ and $\C_{13}$  satisfying 
the defining relations \cite{DGTVZ}
\begin{subequations}\label{eq:BI}
\begin{align}
& \{\C_{12},\C_{23}\}=2(-\C_{13}+\C_{1}\C_{3}+\C_{2}\C_{123}) \,, \label{eq:BI1} \\ 
& \{\C_{12},\C_{13}\}=2(-\C_{23}+\C_{2}\C_{3}+\C_{1}\C_{123}) \,, \label{eq:BI2} \\
& \{\C_{23},\C_{13}\}=2(-\C_{12}+\C_{1}\C_{2}+\C_{3}\C_{123}) \,, \label{eq:BI3}
\end{align}
\end{subequations}
where $\{X,Y\}=XY+YX$.

The algebra $BI(3)$ was first introduced in \cite{TVZ} as an encoding of the bispectral properties of the eponym orthogonal polynomials \cite{BI}. 
In this context the generators $\C_{12}$ and $\C_{23}$ are realized by the Dunkl shift operators of which the Bannai--Ito polynomials are eigenfunctions 
and the operator multiplication by the argument of those polynomials. In this representation the central terms $(\C_1, \C_2, \C_3, \C_{123})$ become 
constants related to the four parameters of the polynomials. 

The centrally extended $BI(3)$ was subsequently defined in \cite{GVZa} following the observation that the Bannai--Ito polynomials are 
essentially the Racah coefficients of the Lie superalgebra $\mathfrak{osp}(1|2)$. This casts $BI(3)$ as the centralizer of the action 
of $\mathfrak{osp}(1|2)$ in the three-fold tensor product $U(\mathfrak{osp}(1|2))^{\otimes 3}$ where $U(\mathfrak{osp}(1|2))$ stands for the 
universal enveloping algebra of $\mathfrak{osp}(1|2)$. The generators $\C_{12}, \C_{23}, \C_{13}$ are then mapped to the Casimir elements attached to 
embeddings of $\mathfrak{osp}(1|2)$ into $\mathfrak{osp}(1|2)^{\otimes 3}$ which are indexed by the 2-element subsets of $\{1, 2, 3\}$. This paved the 
way to the construction of the extension $BI(n)$ of arbitrary rank as the centralizer of the action of $ \mathfrak{osp}(1|2)$ in the $n$-fold 
product $U(\mathfrak{osp}(1|2))^{\otimes n}$ with the generators identified as the Casimir elements associated to $\mathfrak{osp}(1|2)$ 
embeddings now labeled by subsets $A$ of $[n] = \{1, 2, ..., n\}$. This was actually achieved using models of $\mathfrak{osp}(1|2)$ given in terms of Dirac-Dunkl operators \cite{DGV2}, \cite{DGV}.
For reviews of these algebras and some of their applications see \cite{DGTVZ}, \cite{DGVV}.

A notable feature of these tensorial constructs is the fact that the embeddings involved do not all correspond to the simple ones where only the 
factors of the $n$-fold product that are enumerated by the elements the sets $A$ enter non-trivially. The proper Casimir elements are in 
some cases associated to modified embeddings. Sorting this out is addressed here. It will be shown that conjugations of the simple 
embeddings by the universal $R$-matrix will in general be required to ensure that the attached Casimir elements belong to the centralizer.

Throughout this paper, we shall use a presentation of $\mathfrak{osp}(1|2)$ that calls upon a grading involution $P$. This $P$ is group-like under the coproduct
and when this framework is used, it enters in the formula for the universal $R$-matrix. In a separate study \cite{GLV}, expressions for the centralizing elements of $\mathfrak{osp}(1|2)$ have been provided in situations where the grade involution admits a refinement as a product of supplementary involutions. This is manifestly the case under embeddings in tensor products. The centralizing 
elements thus given have been shown in \cite{GLV} to coincide with the Casimirs of the modified embeddings. That this should be so will be made clear in the following. 

The description of the Bannai-Ito algebra in the framework of the universal $R$-matrix of $\mathfrak{osp}(1|2)$ has the striking benefit of allowing to fully characterize abstractly  $BI(n)$ for
arbitrary $n$ (in the centralizer view) without recourse to any model. As shall be shown, the centralizing elements associated to subsets $A$ of $[n] = \{1, 2, ..., n\}$ are given through repeated action of
the coproduct on $\mathfrak{osp}(1|2)$ Casimir elements and conjugation by products of braided universal $R$-matrices. With these generators in hand, the structure relations that they verify can be inferred consistently from those of $BI(3)$ (i.e. \eqref{eq:BI1}-\eqref{eq:BI3}) by exploiting properties of the $R$-matrix and of the permutations of the symmetric group  $\mathfrak S_n$.
A definite picture for the generalized Bannai--Ito algebra $BI(n)$ as the centralizer of $\mathfrak{osp}(1|2)$ in 
$U(\mathfrak{osp}(1|2))^{\otimes n}$ is thus obtained. This approach based on the universal $R$-matrix has already contributed  to the understanding of the Askey--Wilson algebra of rank 1
\cite{CVZ} and the advances presented here in the description of the Bannai-Ito algebra for $n > 3$ should show the way towards a complete picture of the higher ranks  Askey-Wilson algebras.

This letter will proceed as follows. Section \ref{sec:ops} will offer a short review of  $\mathfrak{osp}(1|2)$ and will focus on the universal $R$-matrix of this Lie superalgebra.
In Section \ref{sec:cent}, the centralizing elements of  $\mathfrak{osp}(1|2)$ in 
$U(\mathfrak{osp}(1|2))^{\otimes 3}$ will be given in terms of Casimir elements and the universal $R$-matrix will be shown to play a key role. 
The connection between that centralizer and $BI(3)$ will moreover be made. Section \ref{sec:BIn} will extend the results to $n>3$  and derive the algebra
homomorphism $ BI(n) \rightarrow U(\osp(1|2))^{\otimes n}$ making essential use of the universal $R$-matrix formalism. Short concluding remarks will follow in Section \ref{sec:conc}.

\section{Properties of the Lie superalgebra $\mathfrak{osp}(1|2)$\label{sec:ops}} 

\subsection{The Lie superalgebra $\mathfrak{osp}(1|2)$ \label{sec:super}}

The superalgebra $\osp(1|2)$ has two odd generators $F^\pm$ and three even generators $H$, $E^\pm$ satisfying the following 
(anti-)commutation relations \cite{kac}
\begin{alignat}{2}
& [H,E^\pm]=\pm E^\pm \,,             &\quad & [E^+,E^-]=2H\,,\\
& [H,F^\pm]=\pm \frac{1}{2} F^\pm \,, &\quad & \{F^+,F^-\}=\frac{1}{2}H\,,\\
& [E^\pm,F^\mp]=-F^\pm \,,            &\quad & \{F^\pm,F^\pm\}=\pm\frac{1}{2} E^\pm\;.
\end{alignat}
The $\mathbb{Z}_2$-grading of $\mathfrak{osp}(1|2)$ can be encoded by the grading involution $P$ satisfying
\begin{equation}
 [P,E^\pm]=0 \;, \qquad  [P,H]=0 \;, \qquad  \{P,F^\pm\}=0 \qquad \text{and} \qquad  P^2=1\;.
\end{equation}
One defines the central element $C$ of $U(\mathfrak{osp}(1|2))$ by \cite{Pinc,Les}
\begin{equation}
 C=8[F^+,F^-]P+P\;.\label{eq:Q}
\end{equation}
The $U(\mathfrak{osp}(1|2))$ algebra is endowed with a coproduct $\Delta$ defined as the algebra homomorphism 
satisfying
\begin{alignat}{2}
& \Delta(E^\pm)=E^\pm \otimes 1 +1 \otimes E^\pm \,, &\qquad &\Delta(H)=H\otimes 1 + 1\otimes H\,, \\ 
& \Delta(F^\pm)=F^\pm \otimes P +1 \otimes F^\pm \,, &\qquad &\Delta(P)=P\otimes P\;.
\end{alignat}
We recall that this comultiplication is coassociative 
\begin{equation}
 (\Delta \otimes \id)\Delta= (\id \otimes \Delta)\Delta \
 .\label{eq:coas}
\end{equation}

\subsection{The universal $R$-matrix of $\mathfrak{osp}(1|2)$ \label{sec:uR}}

The universal $R$-matrix of $\osp(1|2)$ is given by
\begin{equation}
 \cR=\frac{1}{2}(1\otimes 1 +P\otimes 1 +1 \otimes P -P \otimes P)\label{eq:UR} .
\end{equation}
For $x\in U(\mathfrak{osp}(1|2))$, it satisfies
\begin{eqnarray}
 \Delta(x) \cR  = \cR \Delta^{op}(x), \label{eq:RDef}
\end{eqnarray}
where the opposite comultiplication $\Delta^{op}(x)=x^{(2)} \otimes x^{(1)}$ if  $\Delta(x)=x^{(1)} \otimes x^{(2)}$ in the Sweedler's notation.
Let us note that 
\begin{equation}
 \cR^2=1\otimes 1 \quad , \qquad  \cR_{21}=\cR\ .
\end{equation}
The universal $R$-matrix \eqref{eq:UR} satisfies 
\begin{eqnarray}
 (\id \otimes \Delta)\cR=\cR_{12}\cR_{13}\ \quad \text{and}\qquad (\Delta \otimes \id)\cR=\cR_{23}\cR_{13}\ .\label{eq:idDRab}
\end{eqnarray}
It verifies also the Yang--Baxter equation
\begin{equation}
 \cR_{12}\cR_{13}\cR_{23}=\cR_{23}\cR_{13}\cR_{12}\ . \label{YBE}
\end{equation}
We remark that in the case of $\osp(1|2)$, the universal $R$-matrix satisfies $[\cR_{12},\cR_{13}]=0$. However, we shall not use this property in the following
so as to keep the computations performed in this letter more generic and applicable to situations involving algebras other than $\osp(1|2)$.

\section{The Bannai--Ito algebra as the centralizer of $\mathfrak{osp}(1|2)$ in $U(\mathfrak{osp}(1|2))^{\otimes3}$ \label{sec:cent}}

\subsection{Centralizer of the diagonal action of $\mathfrak{osp}(1|2)$ in $U(\mathfrak{osp}(1|2))^{\otimes3}$}

To identify $BI(3)$ as the  centralizer of $\mathfrak{osp}(1|2)$ in $U(\mathfrak{osp}(1|2))^{\otimes3}$, it is appropriate to first look for the centralizing elements $X \in U(\mathfrak{osp}(1|2))^{\otimes3}$ such that
\begin{equation}
 [(\Delta\otimes \id)\Delta(x), X]=0\qquad \text{for } x\in \mathfrak{osp}(1|2).
 \end{equation}
 It is straightforward to observe that the elements
 \begin{eqnarray}
 &&C_1=C\otimes 1 \otimes 1\quad , \qquad C_2=1 \otimes C\otimes 1 \quad , \qquad C_3= 1 \otimes 1\otimes C\ ,\label{eq:Cr1}\\
 && C_{12}=\Delta(C)\otimes 1\quad , \qquad  C_{23}=1 \otimes \Delta(C)\quad , \qquad  C_{123}=(\Delta \otimes \id)\Delta(C)\ .\label{eq:Cr2}
\end{eqnarray}
will be centralizing. Now let $\Delta(C)=C^{(1)}\otimes C^{(2)}$ in the Sweedler's notation and write
\begin{equation}
\overline C_{13}=C^{(1)}\otimes 1 \otimes C^{(2)}.
\end{equation}
At first glance one might think that $\overline C_{13}$ also belongs to the centralizer. It is the Casimir element corresponding to the simple homomorphism
 \begin{eqnarray}
   \mathfrak{osp}(1|2)& \rightarrow & \mathfrak{osp}(1|2)^{\otimes3} \nonumber\\
  x&\mapsto & x^{(1)}\otimes 1 \otimes x^{(2)} \nonumber
 \end{eqnarray}
 with $\Delta(x)=x^{(1)}\otimes x^{(2)}$. This however is not true and is where the universal $R$-matrix comes in.
 
  \begin{prop}\label{prop:centr1} 
 The element\footnote{In what follows 
 we shall keep using the inverse of $\cR$ even though $\cR^{-1}=\cR$ (for $\mathfrak{osp}(1|2)$) to make clear that conjugations are involved.}
 \begin{equation}
 C_{13} = \cR_{32} ^{-1} \overline C_{13} \cR_{32}\label{eq:C13a}
 \end{equation}
 belongs to the centralizer of $\mathfrak{osp}(1|2)$ in $U(\mathfrak{osp}(1|2))^{\otimes3}$.
 \end{prop}
 
 \proof Since the Casimir element is central, we have for $x\in \mathfrak{osp}(1|2)$,
\begin{equation}
 [( \Delta\otimes \id)\Delta(x),\Delta(C)\otimes 1]=0\ .
\end{equation}
Using the coassociativity of the comultiplication \eqref{eq:coas} and conjugating by $\cR_{23}$, transforms the previous relation into
\begin{equation}
 [( \id \otimes \Delta^{op})\Delta(x),\cR_{23}^{-1}(\Delta(C)\otimes 1)\cR_{23}]=0\ .
\end{equation}
Finally, exchanging the spaces $2$ and $3$, one gets that $C_{13}$ is in the centralizer
\begin{equation}
 [( \id \otimes \Delta)\Delta(x),\cR_{32}^{-1}\overline C_{13} \cR_{32}]= [( \id \otimes \Delta)\Delta(x), C_{23}] =0\ .
\end{equation}
\endproof
Let us emphasize that $C_{13}$ is in the centralizer whereas $\overline C_{13}$ is not. In particular, for $x=C$ in the previous relation, we get
\begin{equation}
	[C_{123},C_{13}]=0 \ .
\end{equation}
There is the following alternative formula for $C_{13}$.

\begin{prop}\label{prop:centr2}
The element $C_{13}$ is also given by
\begin{equation}
C_{13} = \cR_{12}\overline C_{13} \cR_{12}^{-1}.\label{eq:C13b}
\end{equation}
\end{prop}  

\proof   From property \eqref{eq:RDef}, one gets
\begin{eqnarray}
 C_{13}= \cR_{23}^{-1}\cR_{13} \left(C^{(2)}\otimes 1 \otimes C^{(1)}\right) \cR_{13}^{-1}  \cR_{23}\ .
\end{eqnarray}
Using the Yang--Baxter equation \eqref{YBE}, this relation becomes
\begin{equation}
 C_{13}= \cR_{12}\cR_{13} \cR_{23}^{-1} \cR_{12}^{-1} \left(C^{(2)}\otimes 1 \otimes C^{(1)}\right) \cR_{12} \cR_{23}  \cR_{13}^{-1} \cR_{12}^{-1} \ .
\end{equation}
Now, from \eqref{eq:idDRab}, one deduces that $[\Delta(C)\otimes 1, (\Delta\otimes \id)(\cR)]=[\Delta(C)\otimes 1, \cR_{23}\cR_{13}]=0$ and that 
$[\left(C^{(2)}\otimes 1 \otimes C^{(1)}\right), \cR_{12} \cR_{23}]=0$. One then obtains
\begin{equation}
 C_{13}= \cR_{12}\cR_{13} \left(C^{(2)}\otimes 1 \otimes C^{(1)}\right) \cR_{13}^{-1} \cR_{12}^{-1} \ ,
\end{equation}
which after using \eqref{eq:RDef} again yields the desired result.
\endproof
At this point we can introduce two maps $\hat \tau$ and $\check \tau$ from $U(\mathfrak{osp}(1|2))$ to $U(\mathfrak{osp}(1|2))^{\otimes 2}$ by
\begin{equation}
 \hat \tau(x)= \cR^{-1} (1\otimes x) \cR \quad \text{and} \qquad \check \tau(x)= \cR^{-1} (x\otimes 1) \cR \ .\label{eq:tau}
\end{equation}
\begin{corollary}
The following relations hold in $U(\mathfrak{osp}(1|2))^{\otimes 3}$
\begin{equation}
(\id \otimes \hat \tau)\Delta( C)=C_{13} \quad \text{and} \qquad (\check \tau \otimes \id)\Delta(C)=C_{13}
\end{equation}
\end{corollary}
\proof These results follow directly from Propositions \ref{prop:centr1} and \ref{prop:centr2} and the fact that $\cR=\cR^{-1}$.
\endproof
 Using the definitions \eqref{eq:tau} and the universal $R$-matrix \eqref{eq:UR}, one gets
\begin{eqnarray}
&& \hat \tau(P)=1 \otimes P \ , \quad    \hat \tau(F^{\pm})=P \otimes F^{\pm}\ , \quad  \hat \tau(E^{\pm})=I \otimes E^{\pm}\ \quad \hat \tau(H)=I \otimes H\\\
&&  \check \tau(P)=P \otimes 1 \ , \quad    \check \tau(F^{\pm})=F^{\pm} \otimes P \ , \quad \check \tau(E^{\pm})=E^{\pm} \otimes I, \quad \check \tau(H)=H\otimes I.
 \end{eqnarray}
Either more abstractly with the help of eqs. \eqref{eq:idDRab} or using the formulas above, one readily observes that $\hat \tau$ and $\check \tau$ define coactions, that is verify
\begin{equation}
(\id \otimes\hat \tau)\hat \tau = (\Delta \otimes \id)\hat \tau\  \qquad (\check  \tau \otimes \id  ) \check  \tau = (\id \otimes \Delta ) \check  \tau \ . 
\end{equation}
It hence follows that $(\id \otimes \hat \tau)\Delta$ and $(\check \tau \otimes \id)\Delta$ define two different homomorphisms of $U(\mathfrak{osp}(1|2))$ 
into $U(\mathfrak{osp}(1|2))^{\otimes 3}$ which yield for $C$ the same image, namely:
\begin{equation}
 \C_{13}=\left(8[F^+\otimes P \otimes P+1\otimes 1 \otimes F^+,F^-\otimes P \otimes P+1\otimes 1 \otimes F^-]+1\right)P\otimes 1 \otimes P\;.\label{eq:CCC13}
\end{equation}
This can be checked directly by applying both $(\id \otimes \hat \tau)$ and $(\check \tau \otimes \id)$ to
\begin{eqnarray}
&& \Delta(C)=8\left([F^+\otimes P +1 \otimes F^+ , F^-\otimes P +1 \otimes F^-]+1\right)P\otimes P\\
&& =  16 \left(F^-\otimes F^+ - F^+\otimes F^-\right)\ \left(P\otimes 1\right) +8 C\otimes P+P\otimes C-P\otimes P.\label{eq:delC} 
\end{eqnarray}
Note that
\begin{equation}
\hat \tau (C) =1\otimes C \qquad  \text{and} \qquad \check \tau (C) = C\otimes 1.\
\end{equation}
We may hence pick the homomorphism given by $(\check \tau \otimes \id)\Delta$ and identify the three embeddings labelled by the pairs (12), (23) and (13) (see also \cite{GLV}):
\begin{eqnarray}
&& H_{ij} = H_i + H_j, \quad E^{\pm}_{ij} = E^{\pm}_i + E^{\pm}_j,\qquad i,j = 1,2,3,\nonumber \\
&& F^{\pm}_{12} = F^{\pm}_1 P_2 + F^{\pm}_2, \quad F^{\pm}_{23} = F^{\pm}_2 P_3 + F^{\pm}_3, \quad F^{\pm}_{13} = F^{\pm}_1 P_2 P_3 + F^{\pm}_3, \nonumber \\
&& P_{12} = P_1 P_2, \qquad P_{23} = P_2 P_3 ,\qquad P_{13} = P_1 P_3,\
\end{eqnarray}
with the subscripts denoting (as on the $R$-matrix) the factor in the tensor product where the element appears.
The centralizing elements $C_{ij}$ are then simply the Casimir element given by
\begin{equation}
C_{(ij)}= \Bigl( 8[F^{+}_{(ij)}, F^{-}_{(ij)}]+1\Bigr) P_{(ij)} \, \label{eq:ij}
\end{equation}
as is manifest in particular from \eqref{eq:CCC13} and we now understand the reasons for the choice of the $(13)$ embedding. In this notation we have
\begin{equation}
C_i= \Bigl(8[F^{+}_i, F^{-}_i]+1 \Bigr) P_i, \qquad i=1,2,3\ \label{eq:i}
\end{equation}
\begin{equation}
\text{and} \quad C_{123}= \Bigl( 8[F^{+}_{123}, F^{-}_{123}]+1 \Bigr) P_1P_2P_3 \ \label{eq:123}
\end{equation}
\begin{equation}
\text{with} \quad F^{\pm}_{123}=F^{\pm}_1P_2P_3+F^{\pm}_2P_3+F^{\pm}_3.\
\end{equation}

\subsection{The image of $BI(3)$ in $U(\mathfrak{osp}(1|2))^{\otimes3}$}

We wish to identify the Bannai--Ito algebra  $BI(3)$ with relations \eqref{eq:BI1}-\eqref{eq:BI3} by mapping its generators $\C$ with one, two and three indices onto the corresponding $C$. 
To that end we need to obtain the relations between the Casimir elements $C$. Using the formulas \eqref{eq:ij}, \eqref{eq:i}, \eqref{eq:123}, relation \eqref{eq:BI1} is readily verified under $\C \rightarrow C$.

Note that $C_{13}$ could have been taken to be defined by \eqref{eq:BI1} assuming that the Bannai--Ito relations will be realized. 
(This is typically the approach.) 
Given that $C_{12}$ and $C_{23}$ are centralizing $\mathfrak{osp}(1|2)$ in $U(\mathfrak{osp}(1|2))^{\otimes 3}$, it then follows that $C_{13}$ must also be in the centralizer. 
We have here adopted the view point of first identifying the centralizing elements and hence of first defining $C_{13}$, before obtaining the relations between the generators of the centralizers.
Since the tensorial embedding is so far the only approach that has been designed to obtain the higher rank generalization of the Bannai--Ito algebra, having these definitions of the centralizing 
elements proves essential in this respect.

Given the definitions of $C_{12}$, $C_{23}$ and $C_{13}$, as already said, one directly checks that \eqref{eq:BI1} is satisfied. It is then seen, remarkably, that the remaining defining relations of the Bannai--Ito algebra are implied. One has
\begin{equation}
\{C_{12},C_{23}\}=2(-\cR _{12} \overline C_{13} \cR^{-1} _{12}+C_{1}C_{3}+C_{2}C_{123}).
\end{equation}
Interchanging the factors $1$ and $2$ yields
\begin{equation}
\{C_{21},\overline C_{13}\}=2(-\cR _{21} C_{23} \cR ^{-1}_{21}+C_{2}C_{3}+C_{1}C_{213}).
\end{equation}
Mindful that $C_{21}=\Delta ^{op} (C) \otimes 1$ and that $C_{213}=(\Delta ^{op} \otimes 1)\Delta (C)$, upon conjugating with $\cR ^{-1}_{21} = \cR _{12}$, we find
\begin{equation}
\{C_{12}, \cR_{12}\overline C_{13}\cR^{-1}_{12}\}=2(-C_{23} +C_{2}C_{3}+C_{1}C_{123})
\end{equation}
given that $\cR\Delta^{op} = \Delta \cR$. We thus recover \eqref{eq:BI2} from \eqref{eq:BI1}. The defining relation \eqref{eq:BI3} is also obtained from \eqref{eq:BI1} in a similar fashion. In this case one interchanges the factors $2$ and $3$ and makes use of the other expression for $C_{13}$ namely, $C_{13} = \cR^{-1}_{23} \overline C_{13} \cR_{23}$.

In conclusion, given $C_{12}$, $C_{23}$ and once $C_{13}$ has been defined with the help of the universal $R$-matrix, it is a matter of calculation to obtain one relation between these centralizing elements and one sees thereafter that the 
other two defining relations of the Bannai--Ito algebra follow simply from the first one in light of the properties of the generators and their connection to the $R$-matrix.

\section{The higher rank Bannai--Ito algebras \label{sec:BIn}}

In this section, we shall take $n$ be any positive integer and $[n]=\{ 1,2,\dots,n \}$. 
The higher rank universal Bannai--Ito algebra  $BI(n)$  is generated by $\C_A$ for $A\subset [n]$ (by convention $\C_\emptyset=1$) 
and the following defining relations \cite{DGV}, for $A,B\subset [n]$, 
\begin{equation}
 \{\C_{A},\C_{B}\}=2(-\C_{(A\cup B) \setminus (A\cap B)}+\C_{A\setminus (A\cap B)}\C_{B\setminus (A\cap B)}+\C_{A\cap B}\C_{A\cup B}) \,. \label{eq:BIn} \\ 
\end{equation}
Let us remark that there is a factor $(-2)$ between the generators used here and the ones of \cite{DGV} which explains the apparent discrepancy between the defining relations.  
We shall give an image of $BI(n)$ in $U(\mathfrak{osp}(1|2))^{\otimes n}$. For that, we follow the same logic as before and study the centralizer 
of $\mathfrak{osp}(1|2)$ in $U(\mathfrak{osp}(1|2))^{\otimes n}$.

We define by induction $\Delta^{(k)}=(\id \otimes \Delta^{(k-1)}) \Delta$ with $\Delta^{(0)}=\id$ which allows to define, for $1\leq k \leq \ell \leq n$,
\begin{equation}
 C_{k,k+1,\dots \ell}=1^{\otimes (k-1)} \otimes \Delta^{(\ell-k)}(C) \otimes 1^{\otimes (n-\ell)} \ .\label{eq:Cdefs}
\end{equation}
These elements commute with $\Delta^{(n-1)}(x)$ for $x\in \mathfrak{osp}(1|2)$. We thus obtain elements of the centralizer associated to each subset $K\subset [n]$ with successive integers.
We want to find centralizing elements associated to each subset $A\subset [n]$ without restriction.
Let $\mathfrak S_n$ be the permutation group of $n$ objects generated by the transpositions $s_1,s_2,\dots,s_{n-1}$. 
For $s=s_{i_1}s_{i_2}\dots s_{i_p}$ some permutation of $\mathfrak S_n$ (we recall that any permutation can be written as a product of transpositions),
we define the action $\gamma_s$ on $X\in U(\mathfrak{osp}(1|2))^{\otimes n}$ by
\begin{equation}
 \gamma_s(X)=\check \cR_{i_1}\check \cR_{i_2}\dots \check \cR_{i_p} X (\check \cR_{i_1}\check \cR_{i_2}\dots \check \cR_{i_p})^{-1}, \label{eq:ga}
\end{equation}
where
\begin{equation}
 \check \cR_i= \cR_{i,i+1} \sigma_{i,i+1} 
\end{equation}
and $\sigma_{i,i+1}(x_1\otimes\dots\otimes x_i\otimes x_{i+1}\otimes \dots \otimes x_n)=(x_1\otimes\dots \otimes x_{i+1}\otimes x_{i}\otimes \dots \otimes x_n)\sigma_{i,i+1}$.
Such a $\check \cR_i$ is called braided universal $R$-matrix. It satisfies 
\begin{eqnarray}
 \Delta(x) \check \cR  = \check \cR \Delta(x), \label{eq:RDefb}
\end{eqnarray}
and the braided Yang--Baxter equation
\begin{equation}
 \check \cR_{i} \check \cR_{i+1} \check \cR_{i}=\check \cR_{i+1} \check \cR_{i} \check\cR_{i+1}\ . \label{bYBE}
\end{equation}
Let us emphasize that the definition of $\gamma_s$ does not depend on the choice of the decomposition of the permutation $s$ in terms of the transpositions
since the $\check R_i$ and the $s_i$ satisfy the same algebra.

We define the intermediate Casimir element associated to any subset $A\subset [n]$ as follows 
\begin{equation}
 C_A=C_{s(K)}=\gamma_s (C_K) \label{eq:Cdefq}
\end{equation}
where $K\subset [n]$ with successive integers, $C_K$ is defined by \eqref{eq:Cdefs} and $s \in \mathfrak S_n$ is chosen such that
\begin{equation}
 s(K)=s(\{K_1,K_2,\dots,K_k\})=\{s(K_1),s(K_2),\dots,s(K_k)\}=A.
\end{equation}
We remark that if the permutation $s$ leaves the subset $K$ invariant one gets $\gamma_s (C_K)=C_K$.
It is easy to show using \eqref{eq:RDefb} that $C_A$ is in the centralizer given that $C_K$ is in the centralizer as already proved.

The following example shows that there are different ways to compute $C_A$ depending on the set $K$ we start with.

\begin{example} \label{ex13}
 From $s_1(\{2,3\})=\{1,3\}$ or $s_2(\{1,2\})=\{1,3\}$, the definition \eqref{eq:Cdefq} gives for $C_{13}$
 \begin{eqnarray}
  C_{13}&=&\gamma_{s_1} (C_{23})=\check \cR_1 C_{23} \check \cR_1^{-1}=\cR_{12}\sigma_{12} C_{23}\sigma_{12}^{-1} \cR_{12}^{-1}=\cR_{12} \overline C_{13} \cR_{12}^{-1}\\
  &=&\gamma_{s_2}(C_{12})=\check \cR_2 C_{12} \check \cR_2^{-1}=\cR_{23}\sigma_{23} C_{12}\sigma_{23}^{-1} \cR_{23}^{-1}=\cR_{23} \overline C_{13} \cR_{23}^{-1}\ . 
 \end{eqnarray}
We recover the equivalent expressions \eqref{eq:C13a} or \eqref{eq:C13b} of $C_{13}$ given in the previous section (we recall that $\cR_{12}=\cR_{21}=\cR_{12}^{-1}$). 
\end{example}
To have a well-posed definition of $C_A$, such different paths must lead to the same result.
To confirm that, we must prove that for two subsets $K,L\subset [n]$ of successive integers defined by \eqref{eq:Cdefs} the following relation holds
\begin{equation}
 C_{K}=\gamma_s ( C_{L}) \ , \label{eq:KL}
\end{equation}
where $s(L)=K$. 
It is sufficient to prove \eqref{eq:KL} for the sets $L=\{1,2, \dots,\ell\}$ and $K=\{k+1,\dots, k+\ell \}$ to prove it in general.
The following permutation
\begin{equation}
 s=(s_k s_{k+1} \dots s_{k+\ell-1}) \dots (s_2s_3 \dots s_{\ell+1})(s_1 s_2 \dots s_\ell)
\end{equation}
satisfies $s(L)=K$.
Then, from definition \eqref{eq:ga}, one gets
\begin{eqnarray}
 \gamma_s(C_L)&=&(\check \cR_k \check \cR_{k+1} \dots \check \cR_{k+\ell-1}) \dots (\check \cR_1 \check \cR_2 \dots \check \cR_\ell)
 C_L (\check \cR_\ell\dots \check \cR_{2}\check \cR_1) \dots (\check \cR_{k+\ell-1} \dots \check \cR_{k+1}\check \cR_k )\\
 &=&(\cR_{k,k+1} \cR_{k,k+2} \dots \cR_{k,k+\ell}) \dots ( \cR_{1,k+1} \cR_{1,k+2} \dots \cR_{1,k+\ell})\nonumber \\
 &&\quad C_K (\cR_{1,k+\ell}\dots \cR_{1,k+2} \cR_{1,k+1} ) \dots (\cR_{k,k+\ell} \dots \cR_{k,k+2}\cR_{k,k+1} ) \ . \label{eq:ret}
\end{eqnarray}
The last relation has been obtained using the definition of $\check \cR$ and the properties of $\sigma_{i,i+1}$.
Then, noticing that from relation \eqref{eq:idDRab} one gets $(\id^{\otimes k} \otimes \Delta^{(\ell-1)})(\cR_{i,k+1})=   \cR_{i,k+1} \cR_{i,k+2} \dots \cR_{i,k+\ell}$ (for $1\leq i \leq k$) and $(\id^{\otimes k} \otimes \Delta^{(\ell-1)})(C_{k+1})=C_K$,
one obtains $[\cR_{i,k+1} \cR_{i,k+2} \dots \cR_{i,k+\ell},C_K]=0$ which proves \eqref{eq:KL} in view of \eqref{eq:ret}.\\

We are now ready to present the main result of this section.

\begin{prop}
 The map
 \begin{eqnarray}
  BI(n) &\rightarrow& U(\osp(1|2))^{\otimes n}\\
  \C_A & \mapsto &C_A \nonumber
 \end{eqnarray}
is an algebra homomorphism.
\end{prop}
\proof We must prove that the centralizing elements $C_A$ satisfy the relations \eqref{eq:BIn}.
We know from the previous section that one has
\begin{equation}
\{C_{12},C_{23}\}=2(-\cR _{23} \overline C_{13} \cR^{-1}_{23}+C_{1}C_{3}+C_{2}C_{123}).
\end{equation}
By acting with the coproduct successively on this relation, one gets
\begin{equation}
\{C_{KL},C_{LM}\}=2(- \gamma_s(C_{K,k+1,k+2,\dots, k+m}) +C_{K}C_{M}+C_{L}C_{KLM})
\end{equation}
where $s(K,k+1,k+2,\dots, k+m)=KM$ and $K=\{1,\dots k\}$, $L=\{k+1,\dots k+\ell\}$ and $M=\{k+\ell+1,\dots k+\ell+m\}$.
This proves that for the  sets $K$, $L$ and $M$ given above, the $BI(n)$ relations \eqref{eq:BIn} are satisfied by the corresponding centralizing elements. 
We can similarly see relation \eqref{eq:BIn} to hold when $K$, $L$ or $M$ are chosen empty.
Let $s\in \mathfrak S_n$. Using the definition \eqref{eq:ga}, one gets
\begin{equation}
 \gamma_s ( X X')=\gamma_s ( X) \gamma_s( X').
\end{equation}
Then, we have
\begin{equation}
\{C_{s(KL)},C_{s(LM)}\}=2(- C_{s(KM)} +C_{s(K)}C_{s(M)}+C_{s(L)}C_{s(KLM)}).
\end{equation}
We conclude the proof by remarking that $s(KM)=(s(KL) \cup s(LM)) \setminus (s(KL) \cap s(LM))$, $s(K)=s(KL) \setminus  (s(KL) \cap s(LM))$,
$s(M)=s(LM) \setminus  (s(KL) \cap s(LM))$, $s(L)=s(KL) \cap s(LM)$ and $s(KLM)=s(KL) \cup s(LM)$ and by noting that there exist $K$, $L$ and $M$ and $s$ 
such that $s(KL)=A$ and $s(LM)=B$ for any $A,B\subset [n]$. 
\endproof

\section{Conclusions \label{sec:conc}}
This letter has offered a complete description of the Bannai-Ito algebras as centralizers of the diagonal action of $\mathfrak{osp}(1|2)$ in 
$U(\mathfrak{osp}(1|2))^{\otimes n}$ by bringing the universal $R$-matrix to bear on the topic. This has proved most appropriate. In addition
to the elegance it confers to the presentation, this approach gave answers to questions that had so far been unsettled. It provided an intrinsic algebraic definition
of all centralizing elements independently of the defining relations. It also shed light on the specific form of the intermediate 
embeddings of $\mathfrak{osp}(1|2)$ in $\mathfrak{osp}(1|2)^{\otimes n}$ that yield the generators through the associated Casimir elements. 
Importantly, it has entailed a simple constructive derivation of the structure relations of $BI(n)$ satisfied by these generators through bootstrapping
from the relations of $BI(3)$. Another possible merit is that casting Bannai-Ito algebras in this framework might bring experts familiar with universal $R$-matrices
to contribute further to the field and its applications.

This universal $R$-matrix approach has already been applied to the study of the Askey-Wilson algebra $AW(3)$ \cite{GZ} as 
the centralizer of the diagonal  action of $U_q(\mathfrak{sl}(2))$ into its threefold product and has also been seen to hold promises for advancing 
the understanding of the higher rank $AW(n)$ where one is looking at the centralizer of  $U_q(\mathfrak{sl}(2))$ in 
 $U_q(\mathfrak{sl}(2))^{\otimes n}$\cite{CVZ}. While advances have been made on this last front \cite{DV}, \cite{PW}, a complete description of $AW(n)$ is still lacking. We trust that the treatment
 given here of the Bannai-Ito algebra $BI(n)$ using the universal $R$ - matrix might hold the clues towards bringing this quest to a satisfactory conclusion. We hope to report on this in the near future. \\

\bigskip
\noindent  {\bf Acknowledgements:}
We have much benefited from discussions with L.~Frappat, J.~ Gaboriaud and E.~Ragoucy.
N.~Cramp\'e is gratefully holding a CRM--Simons professorship.
The research of L.~Vinet is supported
in part by a Discovery Grant from the Natural Science and Engineering
Research Council (NSERC) of Canada. M.~Zaimi holds a NSERC graduate scholarship.

\end{document}